\documentclass{amsart}

\usepackage{fancyhdr}
\usepackage{lastpage}
\usepackage{stmaryrd,yhmath}

\newcommand{\bdis}{\begin{displaymath}}
\newcommand{\edis}{\end{displaymath}}
\newcommand{\be}{\begin{equation}}
\newcommand{\ee}{\end{equation}}

\newcommand{\mcal}{\mathcal}
\newcommand{\pd}{\partial}

\newtheorem{lemma}[]{Lemma}

\theoremstyle{definition}

\newtheorem{cor}[]{Corollary}

\theoremstyle{remark}

\newtheorem*{mydef1}{{\bf Theorem}}

\numberwithin{equation}{section}

\begin{document}

\title{Riemann hypothesis and some new asymptotically multiplicative integrals which contain the remainder of the prime-counting function $\pi(x)$}

\author{Jan Moser}

\address{Department of Mathematical Analysis and Numerical Mathematics, Comenius University, Mlynska Dolina M105, 842 48 Bratislava, SLOVAKIA}

\email{jan.mozer@fmph.uniba.sk}


\begin{abstract}
A new parametric integral is obtained as a consequence of the Riemann hypothesis. An asymptotic multiplicability is the main property of this
integral.
\end{abstract}

\maketitle

\section{The result}

\subsection{}

Let us remind that

\be\label{1.1}
\pi(x)=\int_0^x\frac{{\rm d}t}{\ln t}+P(x)=\text{li}(x)+P(x) .
\ee

The following theorem holds true.

\begin{mydef1}
On the Riemann hypothesis
\be \label{1.2}
\int_2^\infty\frac{\ln(xe^{-2})}{x^{3/2+\delta}}P(x){\rm d}x=-\frac 1\delta +\mcal{O}(1),\ \delta\in (0,\Delta)
\ee
where $\Delta$ is a sufficiently small fixed value and the $\mcal{O}(1)$-function is bounded on $[0,\Delta]$.
\end{mydef1}

\subsection{}

Let
\be \label{1.3}
\Omega(\delta)=\int_2^\infty\frac{\ln(xe^{-2})}{x^{3/2+\delta}}\{ -P(x)\}{\rm d}x.
\ee
Then we obtain from (\ref{1.2})

\begin{cor}
\be \label{1.4}
\Omega\left(\prod_{k=1}^n\delta_k\right)\sim\prod_{k=1}^n\Omega(\delta_k),\ \delta_k\to 0,\ k=1,\dots ,n ,
\ee
especially,
\bdis
\Omega(\delta)\sim\sqrt[n]{\Omega(\delta^n)} ,
\edis
i.e., the function $\Omega(\delta)$ possesses the property of the asymptotic multiplicability.
\end{cor}

Let
\bdis
n=p_1^{\alpha_1} p_2^{\alpha_2}\dots p_k^{\alpha_k},\ p_1,\dots ,p_k\to\infty
\edis
be the factorization of a natural number $n$. Then we have (see (\ref{1.4})
\bdis
\frac 1n=\prod_{l=1}^k\left(\frac{1}{p_l}\right)^{\alpha_l} \ \Rightarrow \ \Omega\left(\frac 1n\right)\sim \prod_{l=1}^k
\left\{\Omega\left(\frac{1}{p_l}\right)\right\}^{\alpha_l} ,
\edis
for example.

\section{The main formula}

\subsection{}

Let us remind that
\be \label{2.1}
\begin{split}
& \pi(x)=\int_2^{x}\frac{{\rm d}t}{\ln t}+\bar{P}(x)=\text{Li}(x)+\bar{P}(x) , \\
& \bar{P}(x)=P(x)+\text{V.p.}\int_0^2\frac{{\rm d}t}{\ln t}=P(x)+K;\ K\approx 1.04 ,
\end{split}
\ee
(see (\ref{1.1}), (\ref{2.1}) and \cite{1}, p. 3). The following lemma holds true.

\begin{lemma}
\be \label{2.2}
\int_2^\infty \frac{\pd }{\pd x}\{\ln(1-x^{-s})\}\bar{P}(x){\rm d}x=\ln\zeta(s)+\int_2^\infty\frac{\ln(1-x^s)}{\ln x}{\rm d}x ,
\ee
for $\sigma>1,\ s=\sigma+it$.
\begin{proof}
We apply the formula (\cite{2}, p. 2)
\be \label{2.3}
\ln\zeta(s)=s\int_2^x\frac{\pi(x)}{x(x^s-1)}{\rm d}x,\ \sigma>1 .
\ee
Since
\be \label{2.4}
\frac{s}{x(x^s-1)}=\frac{\pd }{\pd x}\{\ln(1-x^{-s})\} ,
\ee
then we obtain from (\ref{2.3})
\be \label{2.5}
\ln\zeta(s)=\int_2^\infty\frac{\pd }{\pd x}\{\ln(1-x^{-s})\}\pi(x){\rm d}x , \ \sigma>1 .
\ee
Next
\bdis
\frac{{\rm d}\text{Li}(x)}{{\rm d}x}=\frac{1}{\ln x},\ \text{Li}(2)=0 ,
\edis
(see (\ref{2.1})) and
\be\label{2.6}
\begin{split}
& \int_2^\infty\frac{\pd}{\pd x}\{\ln(1-x^{-s})\}\text{Li}(x){\rm d}x=\left.\text{Li}(x)\ln(1-x^{-s})\right|_{x=2}^{x=\infty}- \\
& -\int_2^\infty\frac{\ln(1-x^{-s})}{\ln x}{\rm d}x=-\int_2^\infty\frac{\ln(1-x^{-s})}{\ln x}{\rm d}x ,
\end{split}
\ee
then from (\ref{2.5}) by (\ref{2.1}), (\ref{2.6}) the formula (\ref{2.2}) follows.
\end{proof}
\end{lemma}

\section{The differentiation of the formula (2.2)}

Let $\prod_1$ denote the open rectangle generated by the points $1+\delta\pm i,\ \frac{3}{2}\pm i$ ($0<\delta$ is the sufficiently small fixed value).
The following lemma holds true.

\begin{lemma}
\be \label{3.1}
-\int_2^\infty F(x;s)\bar{P}(x){\rm d}x=\frac{\zeta'(s)}{\zeta(s)}+G(s),\ s\in\Pi_1 ,
\ee
where
\be \label{3.2}
\begin{split}
& F(x;s)=\sum_{n=0}^\infty \{ s(n+1)\ln x -1\}x^{-(n+1)s-1} , \\
& G(s)=\sum_{n=0}^\infty\frac{1}{(n+1)s-1}\frac{1}{2^{(n+1)s-1}} .
\end{split}
\ee
\begin{proof}
We use the formula (see (\ref{2.2}), (\ref{2.4}))
\be \label{3.3}
\int_2^\infty\frac{s}{x^s-1}\frac{\bar{P}(x)}{x}{\rm d}x=\sum_{n=1}^\infty\int_{p_n}^{p_{n+1}}\frac{s}{x^s-1}\frac{\bar{P}(x)}{x}{\rm d}x=
\sum_{n=1}^\infty w_n(s) .
\ee
Since $\bar{P}(x)=\mcal{O}(x)$, and $\bar{P}(x),\ x\in (p_n,p_{n+1})$ is continuous ($p_n$ is a prime number, $p_1=2$) then $w_n(s),\ s\in\Pi_1$ is
an analytic function. Next, the uniform convergence of the series in the set $\Pi_1$ follows from the uniform convergence of the integral  (see
(\ref{3.1}) in the set $\Pi_1$. Thus, by the theorem of Weierstrass, the integral in (\ref{3.3}) is an analytic function in $\Pi_1$. Since
\be \label{3.4}
\begin{split}
& \frac{{\rm d}}{{\rm d}s}\frac{s}{x^s-1}=\frac{1}{x^s-1}-\frac{sx^s\ln x}{(x^s-1)^2}=\frac{x^{-s}}{1-x^{-s}}-\frac{sx^{-s}\ln x}{(1-x^{-s})^2}= \\
& -\sum_{n=0}^\infty\{ s(n+1)\ln x - 1\}x^{-(n+1)s} ,
\end{split}
\ee
then we have (see (\ref{3.3}), (\ref{3.4}))
\be \label{3.5}
\begin{split}
& \frac{{\rm d}}{{\rm d}s}\int_2^\infty\frac{s}{x^s-1}\frac{\bar{P}(x)}{x}{\rm d}x=-\int_2^\infty F(x;s)\bar{P}(x){\rm d}x,\ x\in\Pi_1, \\
& F(x;s)=\sum_{n=0}^\infty \{s(n+1)\ln x -1\}x^{-(n+1)s-1} .
\end{split}
\ee
Similarly, we obtain
\be \label{3.6}
\begin{split}
& \frac{{\rm d}}{{\rm d}s}\int_2^\infty\frac{\ln(1-x^{-s})}{\ln x}{\rm d}x=\int_2^\infty\frac{x^{-s}}{1-x^{-s}}{\rm d}x=\int_2^\infty
\left\{ \sum_{n=0}^\infty x^{-(n+1)s}\right\}{\rm d}x= \\
& = \sum_{n=0}^\infty \frac{1}{(n+1)s-1}\frac{1}{2^{(n+1)s-1}}, \ s\in \Pi_1 .
\end{split}
\ee
Finally, we obtain the formula (\ref{3.1}) by the differentiation of (\ref{2.2}), (see (\ref{3.5}), (\ref{3.6})).
\end{proof}
\end{lemma}

\section{Cancelation of the singular parts corresponding to the pole at $s=1$. The pole at $s=\frac 12$}

Let $\Pi_2$ denote the open rectangle generated by the points $\frac 12\pm i, \frac 32\pm i$. The following lemma holds true.

\begin{lemma}
\be \label{4.1}
\int_2^\infty F(x;s)\bar{P}(x){\rm d}x=-\frac{1}{2s-1}+g(s),\ s\in\Pi_1
\ee
where $g(s),\ s\in \Pi_2$ is the analytic and bounded function.
\begin{proof}
Let (see (\ref{3.1}))
\be \label{4.2}
H(s)=\frac{\zeta'(s)}{\zeta(s)}+G(s),\ s\in\Pi_1 .
\ee
First of all (see (\ref{3.2}))
\be \label{4.3}
\begin{split}
& G(s)=\frac{1}{s-1}\frac{1}{2^{s-1}}+\frac{1}{2s-1}\frac{1}{2^{2s-1}}+g_1(s), \\
& g_1(s)=\sum_{n=2}^\infty \frac{1}{(n+1)s-1}\frac{1}{2^{(n+1)s-1}}, \ s\in \Pi_1  .
\end{split}
\ee
Since
\bdis
\begin{split}
& \frac{1}{s-1}\frac{1}{2^{s-1}}=\frac{1}{s-1}e^{-(s-1)\ln 2}=\frac{1}{s-1}\left\{ 1-(s-1)\ln 2+\mcal{O}(|s-1|^2)\right\}= \\
& =\frac{1}{s-1}-\ln 2+\mcal{O}(|s-1|)=\frac{1}{s-1}+g_2(s) ,
\end{split}
\edis
and similarly
\bdis
\frac{1}{2s-1}\frac{1}{2^{s-1}}=\frac{1}{2s-1}-\ln 2+\mcal{O}(|2s-1|)=\frac{1}{2s-1}+g_3(s)
\edis
then (see (\ref{4.2}), (\ref{4.3}))
\be \label{4.4}
\begin{split}
& H(s)=\frac{\zeta'(s)}{\zeta(s)}+\frac{1}{s-1}+\frac{1}{2s-1}+g_4(s),\ s\in\Pi_1 \\
& g_4(s)=g_1(s)+g_2(s)+g_3(s),
\end{split}
\ee
where $g_4(s), s\in\Pi_2$ is the analytic and bounded function. Next, by the known formula
\be \label{4.5}
\frac{\zeta'(s)}{\zeta(s)}+\frac{1}{s-1}=b-\frac 12\frac{\Gamma'}{\Gamma}\left(\frac s2+1\right)+\sum_{n=1}^\infty
\left(\frac{1}{s-\rho_n}+\frac{1}{\rho_n}\right)=g_5(s)
\ee
where $\zeta(\rho_n)=0$, and $g_5(s),\ s\in\Pi_2$ is the analytic and bounded function. Finally, from (\ref{3.1}) by (\ref{4.2})-(\ref{4.5}) the
asertion of the Lemma 3 follows.
\end{proof}
\end{lemma}

\section{The analytic continuation of the formula (\ref{3.1})}

\subsection{}

The following lemma holds true.

\begin{lemma}
On the Riemann hypothesis
\be \label{5.1}
\int_2^\infty F(x;\sigma)\bar{P}(x){\rm d}x=-\frac{1}{2\sigma-1}+g(\sigma),\ \sigma\in \left(\frac 12,\frac 32\right)
\ee
where $g(\sigma),\ \sigma\in [\frac 12,\frac 32]$ is the bounded function.
\begin{proof}
Let $\Pi_2(\delta)$ denote the open rectangle generated by the points $\frac 12+\delta\pm i; \frac 32\pm i$. We obtain from (\ref{3.2})
\be \label{5.2}
|F(x;s)|=\mcal{O}\left(\frac{\ln x}{x^{\sigma+1}}\right) .
\ee
Next, on the Riemann hypothesis, the following estimate of von Koch
\be \label{5.3}
P(x),\bar{P}(x)=\mcal{O}(\sqrt{x}\ln x)
\ee
holds true (comp. (\ref{2.1})). Then we have (see (\ref{5.2}), (\ref{5.3}))
\bdis
\int_2^\infty |F(x;s)\bar{P}(x)|{\rm d}x=\mcal{O}\left(\int_2^\infty x^{-1-\frac \delta 2}{\rm d}x\right)=\mcal{O}\left(\frac 1\delta\right),
\edis
i.e. the integral
\bdis
\int_2^\infty F(x;s)\bar{P}(x){\rm d}x , \ s\in \Pi_2(\delta)
\edis
is the analytic function. Finally, from the formula (\ref{4.1}), $s\in \Pi_1$, ($g(s),\ s\in \Pi_2$ is the analytic function) by the method of
analytic continuation we obtain the formula
\bdis
\int_2^\infty F(x;s)\bar{P}(x){\rm d}x=-\frac{1}{2s-1}+g(s),\ s\in\Pi_2(\delta) ,
\edis
from which the formula (\ref{5.1}) follows.
\end{proof}
\end{lemma}

\subsection{}

Since (see (\ref{5.2}))
\bdis
\int_2^\infty F(x;\sigma){\rm d}x=\mcal{O}\left(\int_2^\infty x^{-\frac 32+\delta}{\rm d}x\right)=\mcal{O}(1),\ \sigma\in \left[\frac 12,\frac 32\right],
\edis
we obtain, putting $\bar{P}(x)=P(x)+K$ (see (\ref{2.1})) in (\ref{5.1}), the following lemma
\begin{lemma}
On the Riemann hypothesis
\be \label{5.4}
\int_2^\infty F(x;\sigma)P(x){\rm d}x=-\frac{1}{2\sigma-1}+\tilde{g}(\sigma),\ \sigma\in \left(\frac 12, \frac 32\right)
\ee
where $\tilde{g}(\sigma),\ \sigma\in [\frac 12,\frac 32]$ is the bounded function.
\end{lemma}

\section{Proof of the Theorem}

Since (see (\ref{3.2}), comp. (\ref{5.2}), (\ref{5.3}))
\bdis
\begin{split}
& F(x;\sigma)=\frac{\sigma\ln x -1}{x^{\sigma+1}}+\mcal{O}\left(\frac{\ln x}{x^{2\sigma+1}}\right), \\
& \int_2^\infty\frac{\ln x}{x^{2\sigma+1}}|P(x)|{\rm d}x=\mcal{O}\left(\int_2^\infty x^{-\frac 32+\delta}{\rm d}x\right)=\mcal{O}(1),\
\sigma\in \left[\frac 12,\frac 32\right], \\
& \int_2^\infty F(x;\sigma)P(x){\rm d}x=\int_2^\infty\frac{\sigma\ln x-1}{x^{\sigma+1}}P(x){\rm d}x+\mcal{O}(1),\ \sigma\in
\left(\frac 12,\frac 32\right)
\end{split}
\edis
then (see (\ref{5.4}))
\bdis
\int_2^\infty \frac{\sigma\ln x -1}{x^{\sigma+1}}P(x){\rm d}x=-\frac{1}{2\sigma-1}+\mcal{O}(1),\ \sigma\in \left(\frac 12,\frac 32\right) .
\edis
Putting here $\sigma=\frac 12+\delta,\ \delta\in (0,\Delta),\ \Delta<1$, we obtain the formula
\be \label{6.1}
\int_2^\infty\frac{\left(\frac 12+\delta\right)\ln x -1}{x^{\frac 32+\delta}}P(x){\rm d}x=-\frac{1}{\delta}+\mcal{O}(1) .
\ee
However, (see (\ref{5.3}))
\be \label{6.2}
\begin{split}
& \int_2^\infty\frac{\left(\frac 12+\delta\right)\ln x -1}{x^{\frac 32+\delta}}P(x){\rm d}x-
\int_2^\infty\frac{\frac 12\ln x -1}{x^{\frac 32+\delta}}P(x){\rm d}x= \\
& =\delta\int_2^\infty\frac{\ln x}{x^{\frac 32+\delta}}P(x){\rm d}x=\mcal{O}\left(\delta\int_2^\infty\frac{\ln^2 x}{x^{1+\delta}}{\rm d}x\right)= \\
& \mcal{O}\left(\delta\int_2^\infty x^{-1-\frac \delta 2}{\rm d}x\right)=\mcal{O}\left(\delta\frac 1\delta\right)=\mcal{O}(1),\
\delta\in (0,\Delta) .
\end{split}
\ee
The formula (\ref{1.2}) follows from (\ref{6.1}) by (\ref{6.2}).

\thanks{I would like to thank Michal Demetrian for helping me with the electronic version of this work.}

\end{document}